\documentclass[12pt]{article}
\usepackage{amsmath}
\usepackage{amsfonts}
\usepackage{amsthm}
\usepackage{graphicx}
\usepackage{caption}
\usepackage{subcaption}
\usepackage{enumerate}
\usepackage{hyperref}
\usepackage[parfill]{parskip} 

\usepackage{float}
\usepackage{enumitem}

\usepackage{fancyhdr}
\numberwithin{equation}{section}

\DeclareMathOperator{\coeff}{Coeff}
\DeclareMathOperator{\erf}{erf}
\usepackage[labelformat=simple]{subcaption}

\begin{document}

\title{A Guide to the Risk-Averse Gambler and Resolving the St. Petersburg Paradox Once and For All}
\author{Lucy Martinez and Doron Zeilberger}

\maketitle

\begin{abstract}
We use three kinds of computations: Simulation, Numeric, and Symbolic, to guide risk-averse gamblers in general, and offer
particular advice on how to resolve the famous St. Petersburg paradox.

\end{abstract}
\leavevmode
\\
\\
\\
\section{The Famous Saint Petersburg Paradox}

\noindent In the original `infinitarian' version of the St. Petersburg paradox \cite{dM}, a gambler, let's call him Nick,  is tossing a fair coin. If it
lands on Heads, he  gets two ducats, and has to leave the casino. Otherwise he stays, and tosses the coin again, and if it lands on Heads, he gets four ducats,
and has to leave. The reward doubles each time, while he stays at the casino, so if he lasted $k$ rounds, he takes home  $2^k$ ducats.

{\bf Question:} How much should Nick be willing to pay as `entrance fee' to the casino?

Nick's expected gain is
\[\frac{1}{2} \cdot 2 + \frac{1}{4}\cdot 4 + \frac{1}{8} \cdot 8 + \dots \,=\, \sum_{i=1}^{\infty} \frac{1}{2^i}  \cdot 2^i  \, = \, \sum_{i=1}^{\infty} 1 \, = \, \infty.\]

\noindent So Nick should be willing to pay any amount, $M$, (even a billion ducats), since his expected gain would still be $\infty-M=\infty$.

\noindent Obviously, Nick should only be willing to pay a small amount for the privilege of playing. This is the original version of the famous St. Petersburg paradox, that puzzled some of the best minds in probability and economics. Just looking at the references of Wikipedia one can see (in addition to other luminaries, including Laplace), three Nobel Prize winners in economics (Kenneth Arrow, Robert Aumann, and Paul Samuelson).

\subsection{Supporting Maple package and output}

All the results in this article were obtained by the use of the Maple package

\url{https://sites.math.rutgers.edu/~zeilberg/tokhniot/StPete.txt} 

\noindent that also requires the data set

\url{https://sites.math.rutgers.edu/~zeilberg/tokhniot/StPeteData.txt}

\noindent (in the same directory in your computer), whose output files, along with links to diagrams, are available from the front of this article

\url{https://sites.math.rutgers.edu/~zeilberg/mamarim/mamarimhtml/stpete.html}

\subsection{A quick resolution of the St. Petersburg Paradox}

The whole thing is utter nonsense, since it involves an infinite sum, and infinity is a meaningless concept. Besides, life is obviously finite, so the original version of this paradox is just gibberish.

\section{The Finite (and hence Meaningful) Version of the  Saint Petersburg Paradox}

Fix, once and for all, a positive integer $k$, and stipulate that if Nick lasted all the  $k$ rounds, i.e. the coin tosses were all Tails, he would also get $2^k$ ducats, so his
expected gain is
\[
\sum_{i=1}^{k} \frac{1}{2^i}  \cdot 2^i + \frac{1}{2^k} \cdot 2^k \, = \, \sum_{i=1}^{k+1} 1 \, = \, k+1.
\]
Hence, the conventional wisdom of rational choice theory is that he should be willing to pay any amount $n<k+1$ for the privilege of playing, since his expected gain, $k+1-n$ would be positive. Once again, Nick should not accept this bet, if he is only allowed one shot, since his probability of losing money is
very high, and he hates to lose (after all he is {\it risk averse}). Of course, if you want to make money gambling, even if the odds are in your favor, you should
be willing to tolerate some positive chance of losing, but if Nick can insist on being able to play this game many times, then the Central Limit Theorem would guarantee that his chance of exiting the casino a loser can be made as small as he wishes.

{\bf Question}: For a given risk-averseness, i.e. the maximum probability $\epsilon$ of winding up a loser that Nick is willing to take, how many rounds exactly should he insist on?

\section{Simulation}

Stephen Wolfram famously said that formulas and equations are {\it pass\'e}, long live simulation (aka  Monte-Carlo, another casino!). Indeed, in the bad old days, before computers, (poor Count Buffon!), it was impractical to do efficient simulations in real time. In other words, before actually playing for real, have a dry run. Once the gambler decides on insisting that he should be able to repeat the gamble $n$ times, and then he can repeat each such $n$-times game,  $N$ times, and see what happens. The larger the $N$, the better the estimate.

If $n$ is large enough, then he would wind up not losing any money most of the times, but once in a while, he would lose some money, and he hates to lose. He can then count how many of the $N$ `meta-times' were winning, and hence estimate the probability that he will not lose any money with this stipulated $n$.

Now with computers, one can do it very fast, without any calculations! Our Maple package, {\tt StPete.txt} does it for you, dear gambler. First, we have a macro, {\tt StPetePT(n,M)}, that inputs $n$, the number of allowed rounds in one game of the St. Petersburg game, and the ``entrance fee'', $M$. It outputs the probability table of the outcomes of the game. So when $M=0$, it outputs the probability table of outcomes when there is no fee. For the rest of the paper, we remind the reader that Maple syntax requires the user to use a semicolon whenever a macro is used.  For example, trying

{\tt lprint(StPete(6,0));} 

\noindent would output

{\tt [[2, 1/2], [4, 1/4], [8, 1/8], [16, 1/16], [32, 1/32], [32, 1/32]]}

\noindent meaning that with probability $\frac{1}{2}$ you get $2$ ducats, with probability $\frac{1}{4}$ you get $4$ ducats, $\dots$, with probability $\frac{1}{32}$ you get $32$ ducats, and again
with probability $\frac{1}{32}$ you get $32$ ducats (of course, we could have combined these two last outcomes, but for the sake of clarity we prefer to keep it that way).

As noted above, the expected gain is $6$, hence it is still a good deal to have entrance fee $5$. Typing

{\tt lprint(StPete(6,5));} 

\noindent would output

{\tt [[-3, 1/2], [-1, 1/4], [3, 1/8], [11, 1/16], [27, 1/32], [27, 1/32]]}

\noindent meaning that with probability $\frac{1}{2}$ you would lose $3$ ducats, with probability $\frac{1}{4}$ you would lose $1$  ducat, with probability $\frac{1}{8}$ you will win $3$ ducats etc.
Note that if you can only play it once, your probability of losing money is $\frac{3}{4}$, how scary! But with the protection of the law of large numbers, we can try and see what happens, by pure simulation, if you play it many times. Procedure

{\tt Simu1(M,n)} 

\noindent takes any such probability table $M$ and runs the gamble $n$ times and outputs your total gain from this one $n$-times run. Most often, if $n$ is large enough,
you would wind up winning at least some money, but once in awhile you would be a loser. Procedure

{\tt Simu(M,n,N)}

\noindent runs {\tt Simu1(M,n)} $N$ times, and returns the total gain, that should be close to $N$ times the expected gain of $M$ (assumed positive),
followed by the estimated probability that
you will win some money. For example with the above probability table, 

{\tt M:=[[-3, 1/2], [-1, 1/4], [3, 1/8], [11, 1/16], [27, 1/32], [27, 1/32]];}

\noindent typing (once our Maple package {\tt StPete.txt} has been {\tt read}), 

{\tt Simu([[-3, 1/2], [-1, 1/4], [3, 1/8], [11, 1/16], [27, 1/32], [27, 1/32]],100,1000);}

\noindent ouputs something like

{\tt 101.5120000, 0.9150000000}.

\noindent This means that the estimated probability of not losing any money is $0.915$. Of course, this is only an estimate, and every time you get something slightly different. Doing it again, gave us:

{\tt 103.2560000, 0.9200000000}.
 
\noindent We will soon see, using symbolic computation, that the exact probability is $0.9088286275\dots$. So the drawback of simulation (no offense to Wolfram) is that it is only approximate, and also, quite time-consuming, even with a fast computer.

\section{Elementary Symbolic Computation}

Recall that  a probability table for a gamble is a list of pairs of the form:
\[ M=[[M_1,p_1], \dots, [M_r,p_r]].\]

\noindent This means that with probability $p_1$ you would get $M_1$ dollars (or ducats, or whatever),
with probability $p_2$ you would get $M_2$ dollars, $\dots$, and with probability $p_r$ you would get $M_r$ dollars.
In general $M_1, \dots, M_r$ are  any real numbers, but for the sake of simplicity, let's assume that they
are integers. Of course, in real life, currency is discrete, so this assumption is not unrealistic. Note that
some of the $M_i$ are negative, otherwise the decision whether to play would be a no-brainer. If you gamble
you should be willing to lose once in a while. Also, the probabilities $p_i$ are all non-negative, and add-up to one,
\[ p_1 + p_2 + \dots + p_r \, = \, 1.\]

\noindent  The probability generating function (henceforth pgf), in the (formal) variable $x$, is the following Laurent polynomial (i.e. a polynomial that can also have negative exponents, for example $p(x)=1/x+x$),
\[ P_M(x)= \sum_{i=1}^{r} p_i x^{M_i}.\]

\noindent For example, for the above St. Petersburg gamble with six rounds and entrance fee $5$,

{\tt M=[[-3, 1/2], [-1, 1/4], [3, 1/8], [11, 1/16], [27, 1/32], [27, 1/32]]}

\noindent the pgf is
\[ P(x)\,=\, \frac{1}{2} \cdot  x^{-3} \,+\, \frac{1}{4}\cdot x^{-1}+\frac{x^{3}}{8}+\frac{x^{11}}{16}+\frac{x^{27}}{16}. \]

\noindent This is implemented in procedure {\tt PGF(M,x)}, in our Maple package.

Since you are risk-averse,  you are interested in the probability of not losing money, or even better, winning some.
Let's denote by $P(x)^{+}$ the sum of the coefficients whose exponents are positive, then
\[ \left( \sum_{i=1}^{k} p_i x^{M_i} \right )^{(+)} =\sum_{ {{1 \leq i \leq r} \atop {M_i>0}}} p_i. \]

In the above running example, if you only play $M$ once, your probability of winning some money is only $\frac{1}{4}$. But, if you insist on the privilege of playing the gamble a pre-decided number of times, $n$, then your probability
of winning some money is
\[ \left ( P(x)^n \right)^{+}. \]

\noindent  Maple is very good at raising a Laurent polynomial to high powers, expanding them, and then adding up the coefficients of the terms with positive exponents.

\noindent This is implemented in procedure

{\tt ProbPos(M,n)}

\noindent in our Maple package {\tt StPete.txt}. For example, to get the probability of winning some money if you are allowed to gamble $100$ times in the above gamble, type:

{\tt ProbPos([[-3, 1/2], [-1, 1/4], [3, 1/8], [11, 1/16], [27, 1/32], [27, 1/32]],100);}

\noindent getting, immediately that the exact probability is
\[
{\frac{
6125492831448122153753381305179491123116907379470526605886323646825}{6739986666787659948666753771754907668409286105635143120275902562304}}, 
\]
\noindent or more usefully, in decimals,
\[ 0.9088286275, \]
\noindent confirming the estimates that we got above using simulation. So if you play this gamble $100$ times, you know that with probability more than ninety percent you would win some money.
But, being risk-averse, this is too risky! If you insist on being allowed to play exactly $200$ times, then
typing

{\tt evalf(ProbPos([[-3, 1/2], [-1, 1/4], [3, 1/8], [11, 1/16], [27, 1/32], [27, 1/32]],200));}

\noindent would tell you that that the probability of not losing is $0.9733818383$, and if you really want to play it safe, and insist on playing $1000$ times (if you can spare the time)

{\tt evalf(ProbPos([[-3, 1/2], [-1, 1/4], [3, 1/8], [11, 1/16], [27, 1/32], [27, 1/32]],1000));}

\noindent would give you the very reassuring $0.9999947442$.

\subsection{Simplified Gambles}

To simplify matters, and still preserve the St. Petersburg spirit, let's consider the family of gambles whose probability table, let's call it $G_i$, is
\[ G_i:=\left[\left[-1, \frac{i-1}{i}\right], \left[i, \frac{1}{i}\right]\right], \]
whose probability generating function, let's call it $P_i(x)$ is
\[ P_i(x)= \frac{i-1}{i} x^{-1}+ \frac{1}{i} x^i. \]

Note that the expected gain is positive, namely $\frac{1}{i}$, and it is equivalent (changing currency) to the gamble
\[ \left[\left[-i, \frac{i-1}{i}\right], \left[i^2, \frac{1}{i}\right]\right], \]
whose  expected gain is $1$ unit. Let's experiment with $G_{10}$, and the above procedure {\tt ProbPos}.

\noindent If you play $100$ times: the probability of not losing is gotten by typing

{\tt evalf(ProbPos([[-1,9/10],[10,1/10]],100)); }

\noindent giving $0.5487098346$, while if you play $500$ times, typing

 {\tt evalf(ProbPos([[-1,9/10],[10,1/10]],500));}

\noindent would still be the fairly low $0.7453107394$. Wouldn't it be nice if we could compute fast, the sequence \[ \{ProbPos(M,n)\} \quad \]
\noindent for $n \leq N_0$, for any desired $N_0$?

Luckily, thanks to the Almkvist-Zeilberger algorithm \cite{AZ} we can compute many terms very fast, as we will see in the next section.

\section{Advanced Symbolic Computation}

For any Laurent polynomial $P(x)$, we have
\begin{align*}
(P(x))^{(+)} \, &= \,\sum_{j=1}^{\infty} \coeff_{x^j}(P(x))
\, = \,\sum_{j=1}^{\infty} \int_{|x|=1} \frac{1}{2\pi i} \frac{P(x)}{x^{j+1}}\, dx\\
&= \frac{1}{2\pi i} \int_{|x|=1} P(x) \sum_{j=1}^{\infty} \frac{1}{x^{j+1}} \, dx \quad \text{ by Cauchy's residue theorem} \\
&= \frac{1}{2\pi i} \int_{|x|=1} \frac{P(x)}{x(x-1)} \, dx,
\end{align*}
where $\coeff_{x^j}(P(x))$ is the coefficient of $x_j$ in the Laurent polynomial $P(x)$.

Given the one-shot, $M$, with its pgf $P_M(x)$, we are interested in the
probability of winding up with at least some money after $n$ repeats. In other words we are interested in the
sequence
\[ \frac{1}{2\pi i} \int_{|x|=1} \frac{(P_M(x))^n}{x(x-1)} \, dx, \quad n \in \mathbb{N}.\]

Thanks to the Almkvist-Zeilberger algorithm \cite{AZ} (see \cite{D} for a lucid and engaging exposition), such sequences always satisfy a linear recurrence equation with polynomial coefficients. This algorithm is included in {\tt StPete.txt}. The function call is

{\tt OpeProbPos(M,n,Sn)}

\noindent where $M$ is the probability table, $n$ (a symbol!) is the number of repeats, and $Sn$ is the symbol denoting the shift operator in $n$. It also returns the initial conditions. These operators are very complicated, and it is better not to show them to humans. But the computer can use them to compute this sequence very fast. It turns out that eventually, the risk-averse gambler would have to repeat the gamble so many times, that it would be impractical, and he should refuse to play.

\section{Numerics: The Central Limit Theorem to the Rescue}

The advantage of `elementary symbolic computation', and the more efficient and much faster, `advanced symbolic computation' is that
it gives you the exact desired probability. Alas, as the number of repeats $n$ gets larger, sooner or later, it takes too long. It turns out, that for sufficiently large $n$ the approximation given by the Central Limit Theorem gives
you good approximations.

Given a gamble $M=[[M_1,p_1], \dots, [M_r,p_r]]$, define, as usual
\begin{align*}
 \mu &:=\sum_{i=1}^{r} p_i\, M_i, \\
\sigma^2 &:=\sum_{i=1}^{r} p_i\, (M_i-\mu)^2.
\end{align*}

Then the probability that after $n$ repeats, of not losing is approximately
\[ \erf\left(\frac{-\mu\sqrt{n/2}}{\sigma}\right) \]
where $\erf(x)$ is the error-function, built-in in Maple, and $\erf(x)=\frac{2}{\sqrt{\pi}}\int_0^x e^{-t^2} dt$. This is implemented in procedure {\tt ProbPosA}. For small $n$ it is not so good, but it gets better as $n$ gets larger. For example:

{\tt evalf(ProbPos([[-1,9/10],[10,1/10]],100));}

\noindent gives $0.5487098346$, while

{\tt evalf(ProbPosA([[-1,9/10],[10,1/10]],100));}

\noindent gives $0.6190666158$. Not very good!

{\tt evalf(ProbPos([[-1,9/10],[10,1/10]],1000));}

\noindent gives $0.8417618586$, while its Central Limit Theorem approximation gives $0.8310356673$, much better!

{\tt evalf(ProbPos([[-1,9/10],[10,1/10]],10000))}; 

\noindent gives $0.9988718721$, while the approximation is $0.9987784576$, very close! Furthermore the latter is much faster!

\section{Data Files}

Using our Maple package, we prepared lots of useful data files to guide the risk-averse gambler. They are all in the front of this article:

\url{https://sites.math.rutgers.edu/~zeilberg/mamarim/mamarimhtml/stpete.html}.

\noindent It also contains nice pictures. Enjoy!

\section{Graphs}

In this section, we provide graphs for the risk-averseness with different probability tables and gambles.

\begin{figure}[H]
     \centering
     \begin{subfigure}[b]{0.3\textwidth}
         \centering
         \includegraphics[width=\textwidth]{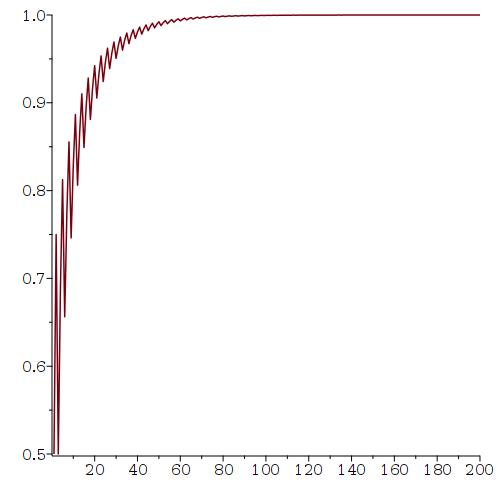}
         \caption{$[[-1,1/2],[2,1/2]]$}
         \label{fig:P2}
     \end{subfigure}
     \hfill
     \begin{subfigure}[b]{0.3\textwidth}
         \centering
         \includegraphics[width=\textwidth]{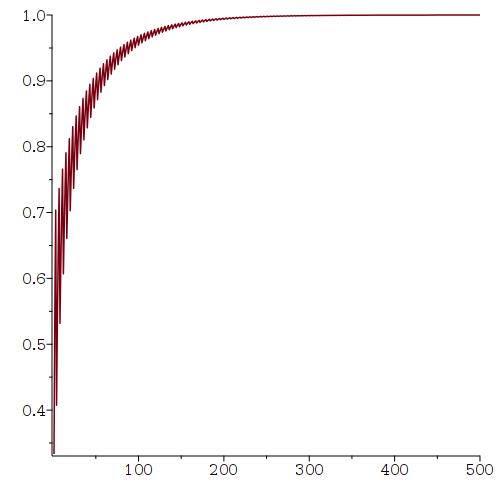}
         \caption{$[[-1,2/3],[3,1/3]]$}
         \label{fig:P3}
     \end{subfigure}
     \hfill
     \begin{subfigure}[b]{0.3\textwidth}
         \centering
         \includegraphics[width=\textwidth]{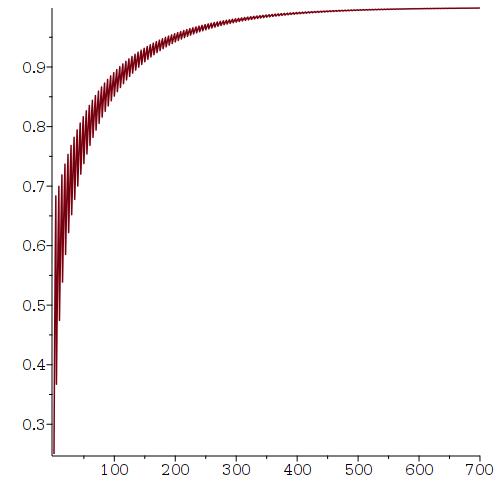}
         \caption{$[[-1,3/4],[4,1/4]]$}
         \label{fig:P4}
     \end{subfigure}
     \hfill
     \begin{subfigure}[b]{0.3\textwidth}
         \centering
         \includegraphics[width=\textwidth]{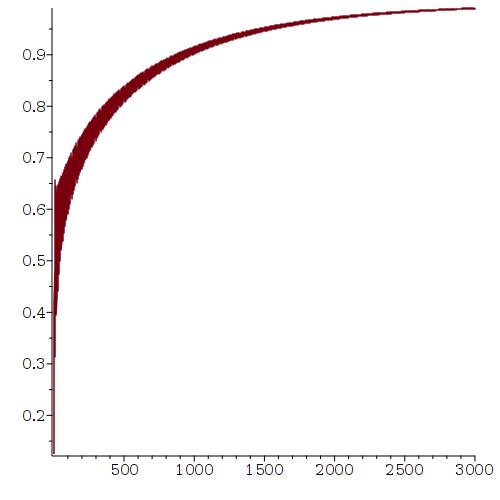}
         \caption{$[[-1,7/8],[8,1/8]]$}
         \label{fig:P8}
     \end{subfigure}
     \hfill
      \begin{subfigure}[b]{0.3\textwidth}
         \centering
         \includegraphics[width=\textwidth]{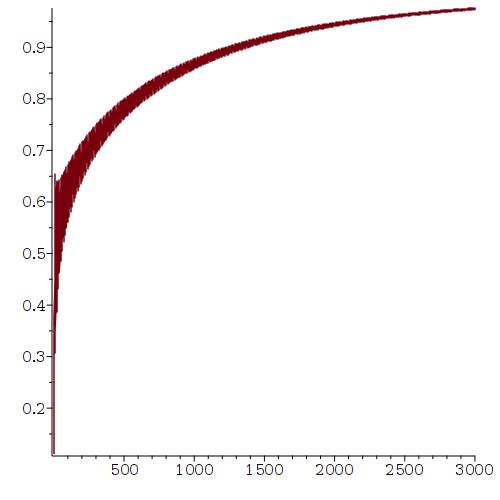}
         \caption{$[[-1,8/9],[9,1/9]]$}
         \label{fig:P9}
     \end{subfigure}
     \hfill
      \begin{subfigure}[b]{0.3\textwidth}
         \centering
         \includegraphics[width=\textwidth]{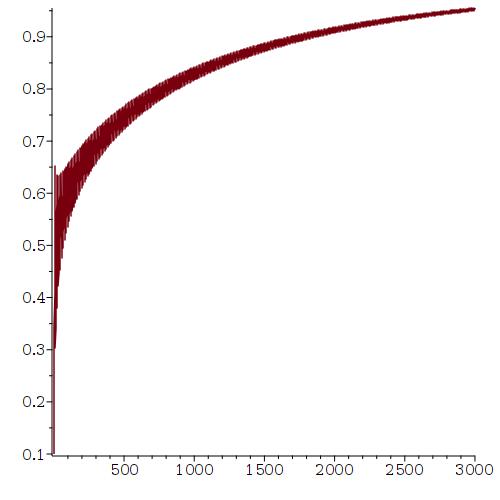}
         \caption{$[[-1,9/10],[10,1/10]]$}
         \label{fig:P10}
     \end{subfigure}
        \caption{The risk-averseness graphs for the corresponding gambles.}
        \label{fig:gamblegraphs}
\end{figure}

\noindent The following list describes each of the graphs in figures \ref{fig:P2}, \ref{fig:P3}, \ref{fig:P4}, \ref{fig:P8}, \ref{fig:P9}, and \ref{fig:P10}:
\begin{enumerate}[leftmargin=*,label=(\alph*)]
\item Given: in one shot, the probability of losing one dollar is $1/2$ and the probability of winning 2 dollars is $1/2$.  The graph illustrates the probability of not losing money if you insist on playing $n$ times, for $n$ from 1 to 200.
\item Given: in one shot, your probability of losing one dollar is $2/3$ and the probability of winning 3 dollars is $1/3$. The graph illustrates the probability of not losing money if you insist on playing $n$ times, for $n$ from 1 to 600.
\item Given: in one shot, your probability of losing one dollar is $3/4$ and the probability of winning 4 dollars is $1/4$. The graph illustrates the probability of not losing money if you insist of playing $n$ times, for $n$ from 1 to 700.
\item Given: in one shot, your probability of losing one dollar is $7/8$ and the probability of winning 8 dollars is $1/8$. The graph illustrates the probability of not losing money if you insist of playing $n$ times, for $n$ from 1 to 3000.
\item Given: in one shot, your probability of losing one dollar is $8/9$ and the probability of winning 9 dollars is $1/9$. The graph illustrates the probability of not losing money if you insist of playing $n$ times, for $n$ from 1 to 3000.
\item Given: in one shot, your probability of losing one dollar is $9/10$ and the probability of winning 10 dollars is $1/10$. The graph illustrates the probability of not losing money if you insist of playing $n$ times, for $n$ from 1 to 3000.
\end{enumerate}

\begin{figure}[H]
     \centering
     \begin{subfigure}[b]{0.3\textwidth}
         \centering
         \includegraphics[width=\textwidth]{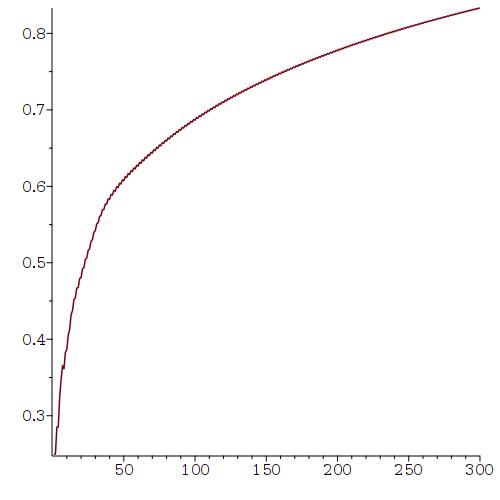}
         \caption{}
         \label{fig:Q7}
     \end{subfigure}
     \hfill
     \begin{subfigure}[b]{0.3\textwidth}
         \centering
         \includegraphics[width=\textwidth]{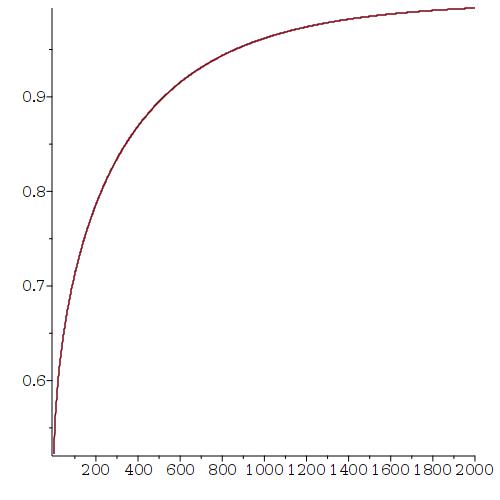}
         \caption{}
         \label{fig:Q7A}
     \end{subfigure}
     \hfill
     \begin{subfigure}[b]{0.3\textwidth}
         \centering
         \includegraphics[width=\textwidth]{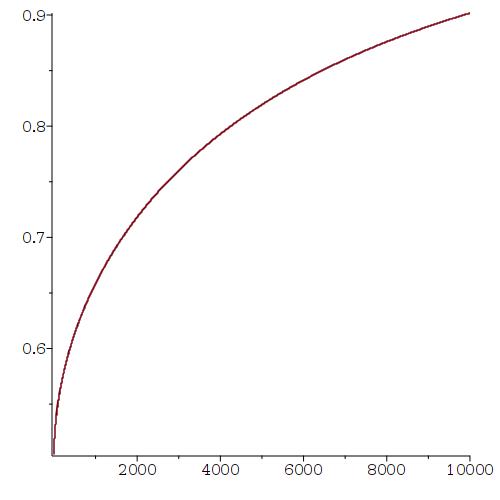}
         \caption{}
         \label{fig:Q11A}
     \end{subfigure}
        \caption{The risk-averseness graphs for the corresponding gambles.}
        \label{fig:centralthmgraphs}
\end{figure}

The following list describes each of the graphs in figures \ref{fig:Q7}, \ref{fig:Q7A} and \ref{fig:Q11A}:
\begin{enumerate}[leftmargin=*,label=(\alph*)]
\item\label{gr:Q7} The risk-averseness graph for the St. Petersburg Gamble with 7 rounds and entrance fee 7. The graph represents the probability of not losing money if you insist on playing $n$ times, for $n$ from 1 to 300 where the probability table is 
\begin{align*}
 &[[-5, 1/2], [-3, 1/4], [1, 1/8], [9, 1/16], [25, 1/32], [57, 1/64], [121, 1/128],\\
 & [121, 1/128]].
\end{align*}
\item\label{gr:Q7A}  The approximate risk-averseness graph for the St. Petersburg Gamble with 7 rounds and entrance fee 7, using the Central Limit Theorem Approximation. The graph represents the probability of not losing money if you insist on playing $n$ times, for $n$ from 1 to 2000 where the probability table is the same as in example \ref{gr:Q7}. 
\item\label{gr:Q11A} The approximate risk-averseness graph for the St. Petersburg Gamble with 7 rounds and entrance fee 7, using the Central Limit Theorem Approximation. The graph represents the probability of not losing money if you insist on playing $n$ times, for $n$ from 1 to 2000 where the probability table is
\begin{align*}
 &[[-9, 1/2], [-7, 1/4], [-3, 1/8], [5, 1/16], [21, 1/32], [53, 1/64], [117, 1/128],\\
 & [245, 1/256], [501, 1/512], [1013, 1/1024], [2037, 1/2048], [2037, 1/2048]].
\end{align*}
\end{enumerate}

\bigskip
\hrule
\bigskip
Lucy Martinez, Department of Mathematics, Rutgers University (New Brunswick), Hill Center-Busch Campus, 110 Frelinghuysen
Rd., Piscataway, NJ 08854-8019, USA. \hfill\break
Email: {\tt  lm1154at math dot rutgers dot edu}    .

Doron Zeilberger, Department of Mathematics, Rutgers University (New Brunswick), Hill Center-Busch Campus, 110 Frelinghuysen
Rd., Piscataway, NJ 08854-8019, USA. \hfill\break
Email: {\tt DoronZeil at gmail dot com}    .

\end{document}